\documentclass[reqno,11pt]{amsart}   	
\usepackage{geometry}                		
\usepackage{graphicx}				
\usepackage{amsmath,amssymb,latexsym,amsfonts}
\usepackage{float}	
\usepackage{cite}									
\allowdisplaybreaks
\begin{document}
\title{Predictive Control of Rabinovich system}
\author[A. Khan, D. Khattar and N. Prajapati]{Ayub Khan, Dinesh Khattar, and Nitish Prajapati*}
\thanks{Received...}
\thanks{A. Khan is with Department of Mathematics, Jamia Millia Islamia, Delhi, India. D. Khattar is with Department of Mathematics, Kirorimal College, University of Delhi, Delhi, India and N. Prajapati is with Department of Mathematics, University of Delhi, Delhi, India.}
\thanks{*Corresponding author. Email: nitishprajapati499@gmail.com}
\subjclass[2000]{34H10, 34H15, 37D45}
\keywords{chaos, rabinovich chaotic system, predictive control, chaos control}

\begin{abstract}
The chaos control problem of continuous time Rabinovich chaotic system is addressed. An instantaneous control input has been designed using predictive control principle to guarantee the convergence of the chaotic trajectory towards an unstable equilibrium point. Numerical simulations are presented to verify all the theoretical analyses. Computational and analytical results are in excellent agreement.
\end{abstract}
\maketitle

\section{Introduction}
Chaos control is one of the major research topic in non linear control systems. Control of chaos refers to the process wherein a tiny perturbation is applied to a chaotic system, in order to realize a desirable behaviour. The first method of control known as OGY method was introduced by Ott, Grebogi and Yorke in 1990 \cite{Ref1}. The main idea was to wait for the natural passage of the chaotic orbit close to the desired periodic behaviour, and then applying a carefully chosen small perturbation, to stabilize such periodic dynamics. Since then the topic has been investigated extensively by many researchers. It has emerged as a new and a very attractive area of research. In this direction many theories and methodologies have been developed, such as time delay feedback method \cite{Ref2}, adaptive control method \cite{Ref3,Ref4}, sliding mode control method \cite{Ref5,Ref6}, occasional proportional feedback method \cite{Ref7}, impulsive control method \cite{Ref8}, backstepping method \cite{Ref9}, passive feedback control method \cite{Ref10}, higher order method \cite{Ref11}, generalized control by non linear high order approach \cite{Ref12}, predictive control method \cite{Ref13,Ref14}.

Recently, predictive control method for controlling continuous chaotic systems has been proposed \cite{Ref14}. Predictive control methods may be considered as a kind of adaptive control strategy \cite{Ref15,Ref16}. The control approach is based on gain feedback configurations. In this method, amplified versions of the difference between future predicted state of uncontrolled chaotic system and the current state is fed back as control input. The stability condition is guaranteed in the same way as in the original OGY method. The predictive control has been applied to the stabilization of various chaotic systems \cite{Ref17,Ref18,Ref19,Ref20,Ref21} and for the synchronization of chaotic satellite systems \cite{Ref22}. The stabilization and synchronization of the chaotic systems have been achieved with great satisfactory performances.

In 1978, Rabinovich differential system was introduced by Pikovski, Rabinovich and Trakhtengerts \cite{Ref23}. This continuous time chaotic system bears a resemblence with the well known Lorenz chaotic system in some properties. The bifurcation diagram and dynamical behaviour of this system has been explored in \cite{Ref24,Ref25}. From Rabinovich system a new 4-D hyperchaotic system was constructed in 2010 \cite{Ref26}. The control of Rabinovich system was first reported in \cite{Ref27} by using controller based on passive control technique. Later in \cite{Ref28}, the control of the system was proposed using sliding mode controller in which control input is required for each state variable.

Motivated by the above chaos control studies, in this paper, we have investigated the control of the continuous time Rabinovich system by using predictive control method. The  advantage of using a predictive controller as an alternative is apparent in the simplicity of configuration and implementation of the control design. The layout of the paper is as follows. Firstly a brief description of predictive control principle is given in Section 2. Then in Section 3 Rabinovich system is described briefly and the controller is designed to stabilize the system using predictive control principle. To verify the theoretical results, numerical simulations have been performed in Section 4 and an excellent agreement between the theoretical and numerical results have been achieved. Finally, concluding remarks are given in Section 5.
\section{Predictive control principle}
Consider an n-dimensional non linear chaotic system described by
\begin{equation} 
\dot{x} (t) = f(x(t)) 
\end{equation}
where $x\in R^n$ is the state vector and we assume that $f$ is differentiable. Our aim is to design a feedback controller $u(t)$ which when added to the dynamical system $(1)$ changes it to the form
\begin{equation} 
\dot{x} (t) = f(x(t))+u(t).
\end{equation}

The controller $u(t)\in R^n$ is designed in such a manner that the trajectory of the controlled system $(2)$ converge to an unstable equilibrium point $x_f$ with only small applied force. The control input $u(t)$ is determined by the difference between the predicted state and the current state as follows
\begin{equation}
u(t) = K(x_p(t)-x(t))
\end{equation}
where K is a gain vector to be determined, $x_p(t)$ is the predicted future state of the uncontrolled chaotic system from the current state $x(t)$. The controlled chaotic system under predictive law is then given by 
\begin{equation}
\dot{x} (t) = f(x(t))+K(x_p(t)-x(t)).
\end{equation}
Using a one-step-ahead-prediction, the predictive controller $(3)$ becomes 
\begin{equation}
u(t) = K(\dot{x} (t)-x(t)). \notag
\end{equation}
The controlled chaotic system $(4)$ then becomes 
\begin{equation}
\dot{x} (t) = (x(t))+K(\dot{x} (t)-x(t)).
\end{equation}

Let $x_f$ be an unstable equilibrium point of system $(1)$. Near $x_f$, we can use the linear approximation for the uncontrolled system by
\begin{equation}
(\dot{x} (t)-x_f) = A(x(t)-x_f)
\end{equation}
where A is the Jacobian matrix of $f(x(t))$ evaluated at the unstable equilibrium point $x_f$, which is defined as follows:
\begin{equation}
A = D_xf(x_f) = {\left.\frac{\partial \dot{x} (t)}{\partial x(t)} \right |}_{x_{f}}. \notag
\end{equation}
Equation $(6)$ is rewritten in the form
\begin{equation}
\delta \dot{x} (t) = A\delta x(t) \notag
\end{equation}
where $\delta x(t) = x(t)-x_f$. The controlled system $(5)$ is linearized around $x_f$ by 
\begin{align}
\delta \dot{x} (t) &=  A \delta x(t)+K(\delta \dot{x} (t)-\delta x(t)) \notag\\
                               &=  A\delta x(t)+K(A\delta x(t)-\delta x(t))\notag\\
                               &= (A+K(A-I))\delta x(t) 
\end{align}
where $I$ is the identity matrix.

Gain K is computed such that $(7)$ is exponentially stable. This implies that the controller gain must satisfy the following inequality
\begin{equation}
|A+K(A-I)|<I. 
\end{equation}  
Gain $K$ exists if and only if det$(A-I)\neq0$. Similar to the OGY method stability is guaranteed in the neighbourhood of equilibrium point $x_f$ \cite{Ref17}. 
\noindent The vicinity of the equilibrium point is determined by 
\begin{equation}
r(t) = |x(t)-x(t-1)|. \notag
\end{equation}
The controlled system is described by 
\begin{align}
\dot{x} (t) &= 
\begin{cases}
  f(x(t))+u(t)    & if \hspace{1mm} r(t) < \epsilon\\
  f(x(t))           & otherwise 
\end{cases}
\end{align}
where $\epsilon$ is a small positive real number.

It is observed that the dynamics of an uncontrolled system performs in such a manner that, given the properties of a chaotic attractor, any trajectory will come, after a possibly long period of time, arbitrarily close to $x_f$. But after the close encounter, the trajectory will rapidly move away from $x_f$. The controller should ensure that once the system is close to the unstable equilibrium point, then it will remain there and asymptotically converge toward $x_f$ \cite{Ref18}.
\section{Predictive control of Rabinovich system}
The Rabinovich chaotic system is defined by the equations \cite{Ref23}
\begin{equation}
\left\{
\begin{aligned}
\dot{x} (t) &=-ax(t)+hy(t)+y(t)z(t),\\ 
\dot{y} (t) &=hx(t)-by(t)-x(t)z(t),\\
\dot{z} (t) &=-dz(t)+x(t)y(t),  
\end{aligned}
\right.
\end{equation}
where $x$, $y$ and $z$ are state variables, and $a$, $b$, $d$ and $h$ are positive real constants. For the parameter values $a=4$, $b=1$, $d=1$ and $h=6.75$, the system $(1)$ displays chaotic behaviour. By simple analysis, according to these parameter values, the equilibrium points of system $(1)$ are: $x_{f_1}=(0,0,0)$, $x_{f_2}=(4.6119,1.3979,6.4469)$ and $x_{f_3}=(-4.6119,-1.3979,6.4469)$. Taking $(1.5, -1.25, 3.5)$ as the initial point, the 3D phase plot of the Rabinovich chaotic system is illustrated in Figure 1. Figure 2 and 3 display the 2D phase plots and the time series respectively.

\begin{figure}[H] 
   \centering
   \includegraphics[height=2.3 in,width=5 in]{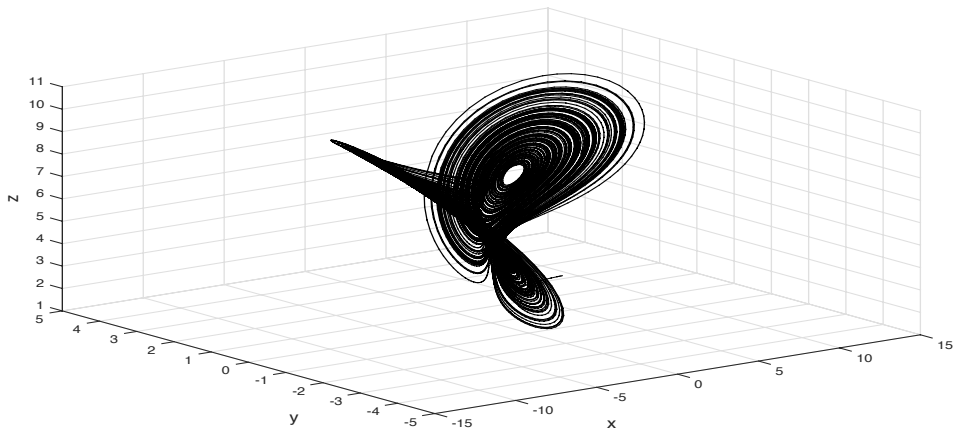} 
   \caption{Phase portrait of Rabinovich chaotic system.}
\end{figure}
\begin{figure}[H] 
   \centering
   \includegraphics[height=8.5 in,width=6.5 in]{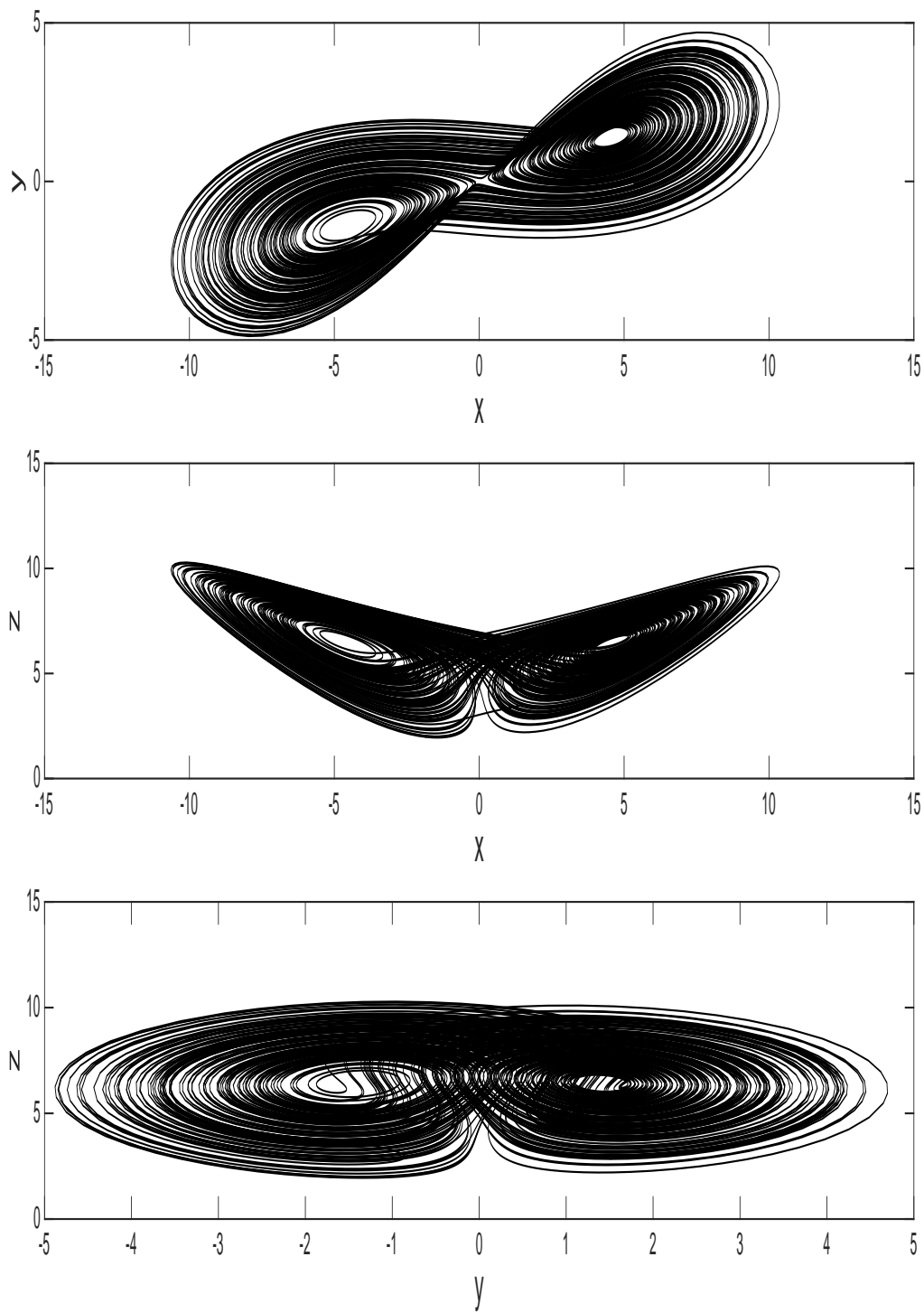} 
   \caption{2D Phase plots of the Rabinovich chaotic system in the $x-y$, $x-z$ and $y-z$ phase planes respectively.}
\end{figure}
\begin{figure}[H] 
   \centering
   \includegraphics[height=4.5 in,width=5.5 in]{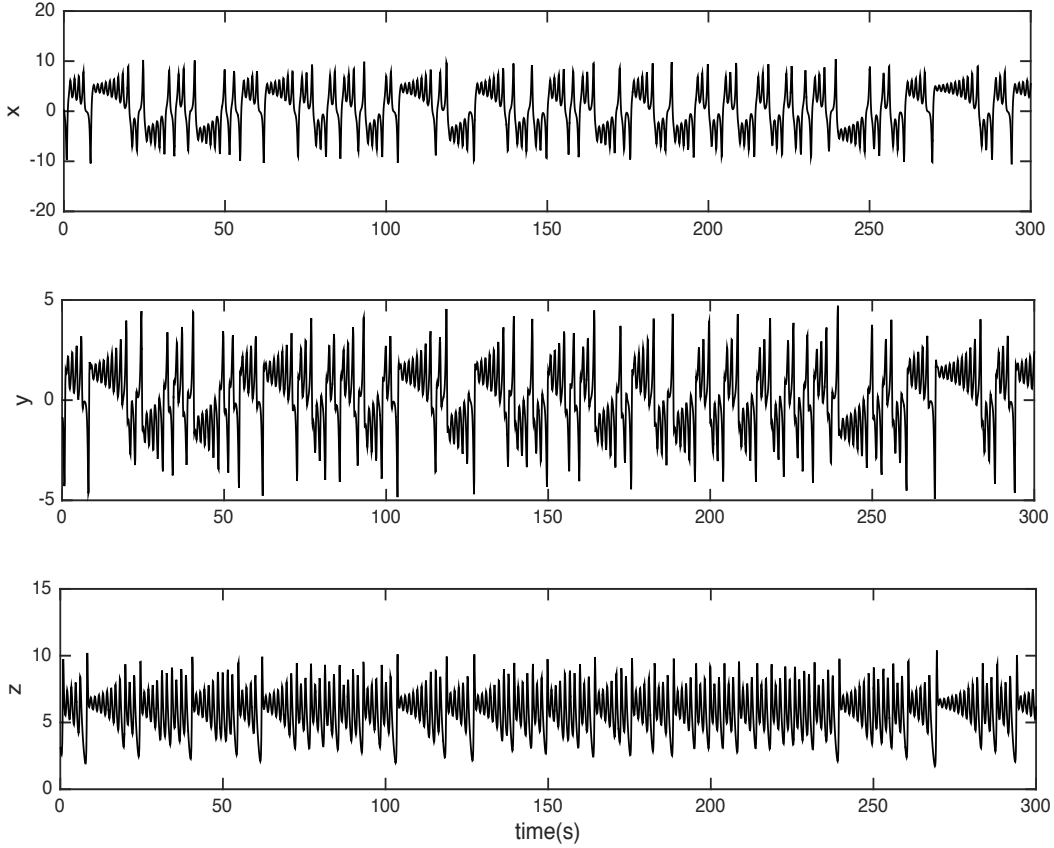} 
   \caption{Time series of the Rabinovich chaotic system for the x, y and z signals.} 
\end{figure}

To stabilize the equilibirium points of Rabinovich system we design a controller based on predictive control principle. After applying the control to the state variable $z$, the system becomes
\begin{equation}
\left\{
\begin{aligned}
\dot{x} (t) &=-ax(t)+hy(t)+y(t)z(t),\\
\dot{y} (t) &=hx(t)-by(t)-x(t)z(t),\\
\dot{z} (t) &=-dz(t)+x(t)y(t)+u(t).
\end{aligned}
\right.
\end{equation}
The controller is designed as follows:
\begin{align}
u(t) &=K(\dot{z} (t)-z(t)) \notag\\
      &=K(-(d+1)z(t)+x(t)y(t)).
\end{align}
where $\dot{z} (t)$ is the predicted future state of the state variable $z$ and $z(t)$ is the current state. Thus, the controlled chaotic system $(11)$ is expressed as 
\begin{equation}
\left\{
\begin{aligned}
\dot{x} (t) &=-ax(t)+hy(t)+y(t)z(t),\\
\dot{y} (t) &=hx(t)-by(t)-x(t)z(t),\\
\dot{z} (t) &=-dz(t)+x(t)y(t)+K(-(d+1)z(t)+x(t)y(t)).
\end{aligned}
\right.
\end{equation}
Linearizing the state variable $z$ around $x_{f_i}$($i = 1, 2, 3$), we obtain
\begin{align}
\delta \dot{z} (t) &={\left.\frac{\partial \dot{z} (t)}{\partial z(t)} \right |}_{x_{f_i}} \delta z(t) \notag \\
                         &=(-d-K(d+1)) \delta z(t).
\end{align}

The controller $(12)$ will stabilize the uncontrolled system $(10)$ in the local neighbourhood of equilibrium point $x_{f_i}$, if there exist a gain $K$ such that $(14)$ is stable. Thus in order to apply the proposed predictive controller we have to determine the gain $K$. For the parameter value $d = 1$ and using the stability criteria as given in $(8)$, the feedback gain $K$ must satisfy the condition
\begin{equation}
|-1-2K|<1. \notag
\end{equation}
Consequently, for $-1<K<0$, the predictive controlled Rabinovich system is obtained with $(13)$.
\section{Numerical Simulations}
We carry out the simulations using the MATLAB. The fourth-order Runge-Kutta method has been used in the algorithm. A time step size 0.1 has been employed to solve the differential equation. For $K=-0.6$, $\epsilon=0.1$ and starting from $(x(0), y(0), z(0))=(1.5, -1.25, 3.5)$, the control input $u$ is applied for $t>40$. After a brief transitory phase, we see, the system trajectory is stabilized around its unstable equilibrium point $x_{f_3}$ as shown in Figure $4$. When control is activated for $t>100$, in this case, the control input $u$ stabilizes the chaotic trajectory on $x_{f_2}$ with extremely small applied force as shown in Figure $5$.
\begin{figure}[H] 
   \centering
   \includegraphics[height= 5.15 in,width=5.5in]{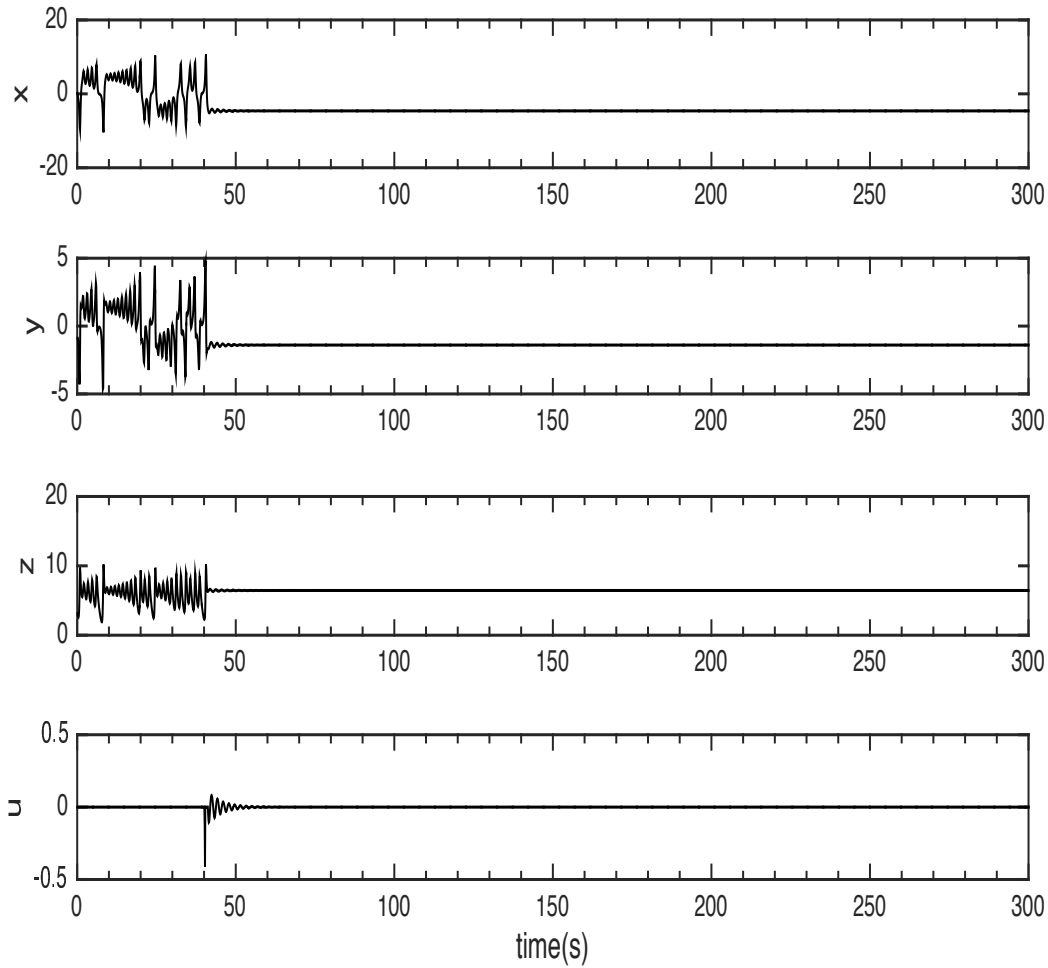}  
   \caption{Stabilization of the Rabinovich chaotic system on the unstable equilibrium point $x_{f_3}$.}
\end{figure}
\begin{figure}[H] 
   \centering
   \includegraphics[height= 5.15 in,width=5.5in]{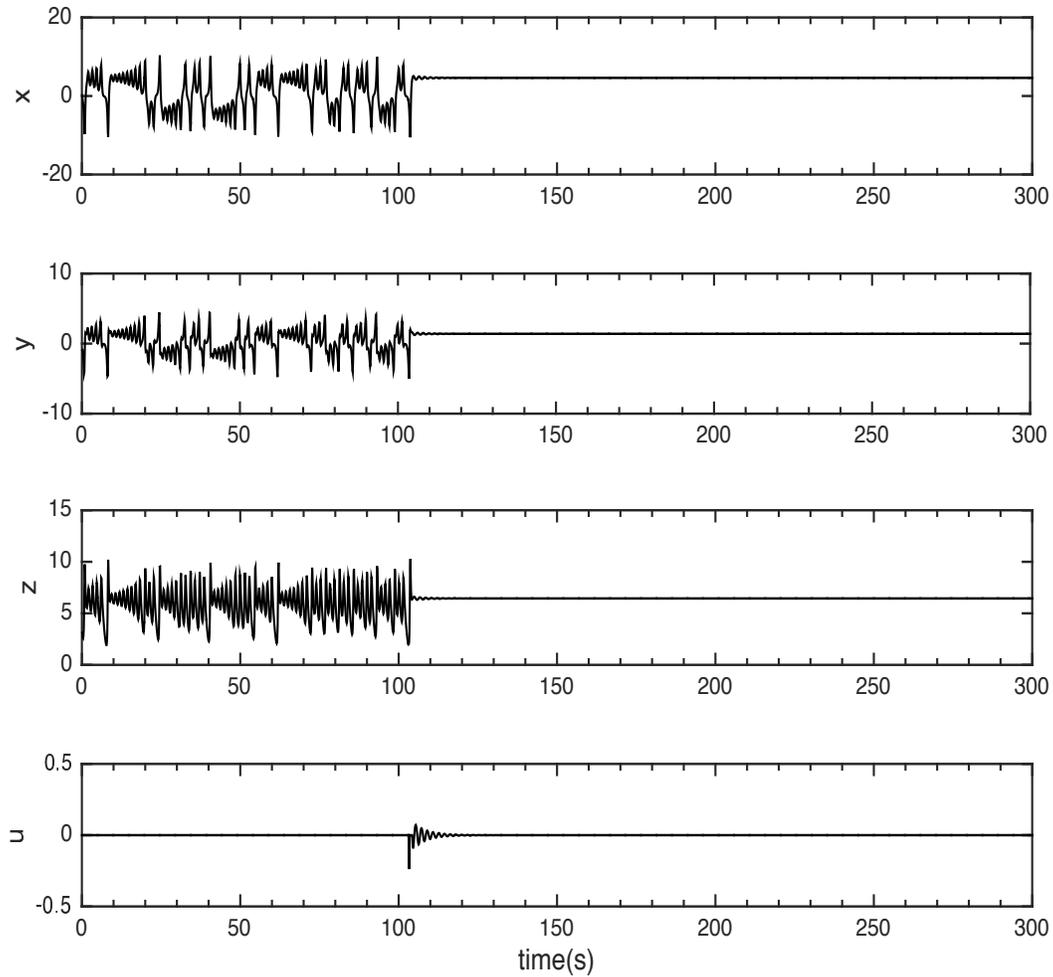}  
   \caption{Stabilization of the Rabinovich chaotic system on the unstable equilibrium point $x_{f_2}$.}
\end{figure}

\section{Conclusion}
This paper investigates the control of the continuous time Rabinovich chaotic system by predictive control method. The controller designed by predictive control method has been used to stabilize the system trajectories on the unstable equilibrium points successfully. The predictive controller designed is simple to configure and implement. It's single controller approach resulting in ease of application and low cost production makes it an attractive alternative to other control methods. Numerical simulations shows the effectiveness of the proposed method. To check the robustness of the predictive method, control of the system with parameter disturbance and uncertainty, will be an interesting future work.
\section*{Acknowledgements}
The work of the third author is supported by the Junior Research Fellowship of Council of Scientific and Industrial Research, India(Grant no. 09/045(1319)/2014-EMR-I).


\begin{thebibliography}{99}
\bibitem{Ref1} E. Ott, C. Grebogi\ and\ J.A. Yorke: Controlling chaos. Physics Review Letters. {\bf 64}, (1990), 1196--1199.
\bibitem{Ref2} K. Pyragas: Continuous control of chaos by self controlling feedback. Physics Letters A. {\bf 170}, (1992), 421--428.
\bibitem{Ref3} M. di Bernardo: An adaptive approach to the control and synchronization of continuous-time chaotic systems. International J. of Bifurcation and Chaos. {\bf 6}, (1996), 557--568.
\bibitem{Ref4} A. Lor\'\i a\ and\ A. Zavala-R\'\i o: Adaptive tracking control of chaotic systems with applications to synchronization. IEEE Transactions on Circuits and Systems I. {\bf 54}, (2007), 2019--2029.
\bibitem{Ref5} J.J. Slotine: Sliding controller design for nonlinear systems. International J. of Control. {\bf 40}, (1984), 421--434.
\bibitem{Ref6} K. K2nishi, M. Hirai\ and\ H. Kokame: Sliding mode control for a class of chaotic systems. Physics Letters A {\bf 245}, (1998), 511--517.
\bibitem{Ref7} N. Inaba\ and\ T. Nitanai: OPF chaos control in a circuit containing a feedback voltage pulse generator. IEEE Transactions on Circuits and Systems I. {\bf 45}, (1998), 473--480.
\bibitem{Ref8} X. Liu: Impulsive stabilization and control of chaotic system. Nonlinear Analysis {\bf 47}, (2001), 1081--1092. 
\bibitem{Ref9} C. Wang\ and\ S.S. Ge: Synchronization of two uncertain chaotic systems via backstepping. International J. of Bifurcation and Chaos. {\bf 11}, (2001), 1743--1751.
\bibitem{Ref10} W. Yu: Passive equivalence of chaos in Lorenz system. IEEE Transactions on Circuits and Systems I. {\bf 46}, (1999), 876--878.
\bibitem{Ref11} A. Boukabou\ and\ N. Mansouri: Controlling chaos in higher-order dynamical systems. International J. of Bifurcation and Chaos. {\bf 14}, (2004), 4019--4025.
\bibitem{Ref12} A. Boukabou\ and\ N. Mekircha: Generalized chaos control and synchronization by nonlinear high-order approach. Mathematics and Computer in Simulation. {\bf 82}, (2012), 2268--2281. 
\bibitem{Ref13} T. Ushio\ and\ S. Yamamoto. Prediction-based control of chaos. Physics Letters A. {\bf 264}, (1999), 30--35.
\bibitem{Ref14} A. Boukabou, A. Chebbah\ and\ N. Mansouri: Predictive control of continuous chaotic systems. International J. of Bifurcation and Chaos. {\bf 18}, (2008), 587--592.
\bibitem{Ref15} R.P. Bithmead, M. Gevers\ and\ V. Wertz: Adaptive optimal control- The thinking GPC. Prentice-Hall, Englewood Cliffs, NJ, 1990.
\bibitem{Ref16} K.S. Park, J.B. Park, Y.H. Choi, T.S. Yoon\ and\ G. Chen: Generalized predictive control of discrete time chaotic systems. International J. of Bifurcation and Chaos. {\bf 8}, (1998), 1591--1597.
\bibitem{Ref17} A. Boukabou\ and\ N. Mansouri: Fuzzy predictive controller for unknown discrete chaotic systems. International J. of Bifurcation and Chaos. {\bf 17}, (2007), 2141--2148.
\bibitem{Ref18} A. Boukabou, B. Sayoud, H. Boumaiza\ and\ N. Mansouri: Control of n-scroll Chua's circuit. International J. of Bifurcation and Chaos. {\bf 19}, (2009), 3813--3822.
\bibitem{Ref19} S. Hadef\ and\ A. Boukabou: Control of multi-scroll Chen system, J. of Franklin Institute. {\bf 351}, (2014), 2728--2741. 
\bibitem{Ref20} D. Sadaoui, A. Boukabou\ and\ S. Hadef: Predictive feedback control and synchronization of hyperchaotic systems. Applied Mathematics and Computation. {\bf 247}, (2014), 235--243.
\bibitem{Ref21} M. Messadi, A. Mellit, K. Kemih\ and\ M. Ghanes: Predictive control of a chaotic permanent magnet synchronous generator in a wind turbine system. Chinese Physics B. {\bf 24}, (2015), 010502.
\bibitem{Ref22} D. Sadaoui, A. Boukabou, N. Merabtine\ and\ M. Benslama: Predictive synchronization of chaotic satellite systems, Expert Systems with Applications. {\bf 38}, (2011), 9041--9045.
\bibitem{Ref23} A.S. Pikovski, M.I. Rabinovich\ and\ V.Y. Trakhtengerts: Onset of stochasticity in decay confinement of parametric stability. Soviet Physics JETP. {\bf 47}, (1978), 715--719.
\bibitem{Ref24} S. Neukirch: Integrals of motion and semipermeable surfaces to bound the amplitude of plasma instability. Physics Review E. {\bf 63}, (2001), 036202.
\bibitem{Ref25} J. Llibre, M. Messias\ and\ P.R. Silva: On the global dynamics of the Rabinovich system. J. of Physics A. {\bf 41}, (2008), 275210.
\bibitem{Ref26} Y. Liu, Q. Yang\ and\ G. Pang: A hyperchaotic system from the Rabinovich system. J. of Computational and Applied Mathematics. {\bf 234}, (2010), 101--113.
\bibitem{Ref27} S. Emiro\u glu\ and\ Y. Uyaro\u glu: Control of Rabinovich chaotic system based on passive control. Scientific Research and Essays. {\bf 5}, (2013), 3298--3305.
\bibitem{Ref28} U. E. Kocamaz, Y. Uyaro\u glu\ and\ H. Kizmaz: Control of Rabinovich chaotic system using sliding mode control. International J. of Adaptive Control and Signal Processing. {\bf 28}, (2014), 1413--1421.

\end{thebibliography}
\end{document}